\documentclass[12pt]{article}
\usepackage{amsmath,amssymb,amsfonts,euscript}
\title{{\large\bf A comment to: On $3$-colorable planar\\ graphs without short cycles }}
\author{ S. Akbari $^{b,a}$ and Behrooz Bagheri Gh.${}^a$}
\date{}

\begin{document}
\maketitle
\begin{center}
$a$ Department of Mathematical Sciences \\
Sharif University of Technology \\
\vspace*{5mm}
$b$ 
 Institute for Studies in Theoretical Physics and Mathematics (IPM)\\

\end{center}

\begin{abstract}
Let $G$ be a graph. It was  proved that if $G$ is a planar graph
without $\{4,6,7\}$-cycles  and without two $5$-cycles sharing
exactly one edge, then $G$ $3$-colorable. We observed that the proof
of this result is not correct.

\end{abstract}


\section{}        
 Let $G$ be a simple  graph with vertex set $G$. A {\it planar} graph is one that can be drawn on a plane in such
a way that there are no ``edge crossings," i.e. edges intersect only
at their common vertices. A $k$-{\it face} is a face whose boundary
has  $k$ edges. A cycle $C$ in a planar graph $G$ is said to be {\it
separating}, if $int(C)\neq \varnothing $ and $ext(C)\neq
\varnothing$,
 where $int(C)$ and $ext(C)$ denote the sets of vertices located
inside and outside $C$, respectively.  Let $C_i$ denote an
$i$-cycle. A $k$-{\it coloring} of $G$ is a mapping $c$ from $V(G)$
to the set $\{1,\ldots,k\}$ such that $c(x)\ne c(y)$ for any
adjacent vertices $x$ and $y$. The graph $G$ is $k$-{\it colorable}
if it has a $k$-coloring. Let $\cal G$ denote the class of planar
graphs without $\{4,6,7\}$-cycles and without two $5$-cycles sharing
exactly one edge.\\

\noindent The following theorem was proved in \cite{che}.\\

\noindent {\bf Theorem.} Let $G$ be a graph in $\cal G$ that
contains $5$-cycles. Then every proper $3$-coloring of the vertices
of any $3$-face or $9$-face of $G$ can be extended into a proper
$3$-coloring of the whole graph.\\

The authors divided the proof of the theorem into $9$ lemmas and
their proof is by contradiction. They considered a graph $G$ with
the minimum number of vertices such that satisfies the assumptions,
but the assertion of the theorem is not true for $G$. In the proof
of many lemmas they made a common mistake. For instance in Lemma 1,
they applied induction on $G \setminus int(C_i)$, where $C_i$ is a
separating $i$-cycle. But $G \setminus int(C_i)$ has no necessarily
$5$-cycle and so their proof is not correct.
{}

\end{document}